\documentclass{article}

\newcommand{\R}{I\!\! R}

\newcommand{\Z}{I\!\! Z}
\newcommand{\N}{I\!\! N}
\newcommand{\Q}{I\!\! Q}

\begin{document}

Aristide Tsemo,

Department of Mathematics and Computer Science,

   Ryerson University     350, Victoria Street

        Toronto, ON

        M5B 2K3

        Canada   taristid@scs.ryerson.ca

\bigskip
\bigskip

 \centerline{\bf G\'eom\'etrie affine, g\'eom\'etrie symplectique.}

\bigskip
\bigskip

\centerline{\bf Abstract.}

\bigskip

{\it In this paper we study the relation between affine manifolds
and symplectic geometry, Liberman and Weinstein have shown that a
leaf of a lagrangian foliation is endowed with an affine
structure, we translate properties, and conjectures of affine
manifolds in symplectic geometry. This allows to show the
Auslander to be true if the linear holonomy is contained in
$Gl(n,{\Z})$.}

\bigskip
\bigskip

Une vari\'et\'e symplectique $(M,\omega)$ de dimension $2n$, est
une vari\'et\'e diff\'erentiable $M$, munie d'une $2-$forme
ferm\'ee $\omega$, telle que $\Lambda^n\omega$ soit une forme
volume de $M$. Un feuilletage ${\cal F}$ de dimension $n$, de  $M$
est lagrangien si et seulement si les feuilles de ${\cal L}$ sont
isotropes pour $\omega$, ce qui est \'equivalent \`a dire que la
restriction de la forme symplectique sur chacune des feuilles est
nulle.

Weinstein a montr\'e que les feuilles d'un feuilletage lagrangien
sont munies d'une structure de vari\'et\'es affines, i.e d'une
connexion dont les tenseurs de courbure et de torsion sont nulles.
Une telle construction similaire avait \'et\'e effectu\'ee par
Paulette Libermann.

Rappelons la construction de Dazord de la connection de
Libermann-Weinstein: Soient $L_0$ une feuille de ${\cal L}$,
$\chi({\cal L})$ l'alg\`ebre des champs de vecteurs tangents \`a
${\cal L}$, et $\chi(n({\cal L}))$ l'ensemble des sections du
fibr\'e normal de ${\cal L}$.

La connexion de Bott du feuilletage ${\cal L}$ est d\'efinie sur
$\chi({\cal L})\otimes\chi(n({\cal L}))$ par:

$$
\hat \nabla_XY=u([X,Y'])
$$
o\`u $Y'$ est un \'el\'ement de $\chi(M)$ au-dessus de $Y$, et $u$
la projection canonique $\chi(M)\rightarrow \chi(n({\cal L}))$.
Cette connexion induit sur le dual $(\chi(n({\cal L})))^*$ une
connection $\nabla'$ d\'efinie par la formule classique:

$$
\nabla'_Xf(t)=L_X(f(t)) - f(\hat\nabla_Xt)
$$

La dualit\'e symplectique permet d'identifier $\chi({\cal L})$ \`a
$(\chi(n({\cal L})))^*$ puisque ${\cal L}$ est un feuilletage
lagrangien. Par cette identification, la connection $\nabla'$
induit sur chaque feuille $L_0$ de ${\cal L}$ une connexion
$\nabla_{L_0}$ dont les tenseurs de courbure et de torsion sont
nuls. Ceci est \'equivalent \`a d\'efinir la structure
diff\'erentielle de $L_0$ par un atlas dont les fonctions de
transitions sont des applications affines.

Consid\'erons un syst\`eme de coordonn\'ees de Darboux
$(q_1,..,q_n,p_1,..,p_n)$ adapt\'e au feuilletage, dans ce
syst\`eme le feuilletage est d\'efini par les \'equations
$dq_1=...=dq_n=0$. Les coordonn\'ees $(p_1,..,p_n)$ d\'efinissent
la structure affine de Libermann-Weinstein.

La connexion $\nabla_{L_0}$ est le dual de la connexion de Bott.
On en d\'eduit ainsi un isomorphisme entre l'holonomie de
$\nabla_{L_0}$ et l'holonomie infinit\'esimale de ${\cal L}$ en
$L_0$. Le but de ce papier et de transcrire les propri\'et\'es des
vari\'et\'es affines en utilisant le dictionnaire d\'efini par la
g\'eom\'etrie lagrangienne. La topologie des vari\'et\'es
symplectiques compactes a \'et\'e \'etudi\'ee par divers auteurs,
Donaldson a montr\'e l'existence de sous-vari\'et\'es
symplectiques en toutes codimensions, Auroux a \'etudi\'e
l'existence de d\'eformations continues de ces vari\'et\'es
symplectiques. On se sert de ces constructions pour d\'emontrer la
conjecture d'Auslander lorsque l'holonomie lin\'eaire est contenue
dans $Gl(n,{\Z})$. Puis nous  d'\'etendons le principe de
r\'eduction de Molino en utilisant la notion de tour de torseurs.

\bigskip

{\bf Conjecture d'Auslander et croissance.}

\medskip

Une vari\'et\'e affine $(M,\nabla_M)$ est compl\`ete si et
seulement si la connexion $\nabla_M$ est g\'eod\'esiquement
compl\`ete. L. Auslander a conjectur\'e que le groupe fondamental
d'une vari\'et\'e affine compacte et compl\`ete $(M,\nabla_M)$ est
polycyclique. Cette conjecture a faite l'objet de nombreux
travaux, et a \'et\'e prouv\'ee dans les cas suivants:

- La dimension de $M$ est inf\'erieure \`a $3$ par Fried et
Goldman,

- La dimension de $M$ est inf\'erieure \`a $6$, par Abels, Soifer
et Margulis.

- L'holonomie de la connection $\nabla_M$ est incluse dans
$O(n-1,1)$, par Goldman et Kamishima.

- Le groupe des transformations affines de $(M,\nabla_M)$ est de
cohomog\'eneit\'e sup\'erieure \`a $1$ par Tsemo.

\medskip

{\bf Definitions 2.1} Soient $(M,{\cal L})$ une vari\'et\'e
diff\'erentiable compacte $M$, munie d'un feuilletage ${\cal L}$,
et $(U_i)_{i\in I}$ un atlas de $M$ tel que la restriction de
${\cal L}$ \`a $U_i$ est simple. Consid\'erons $L_0$ une feuille
de ${\cal L}$, une plaque $P_i$ de $L_0$. est une composante
connexe de $L_0\cap U_i$. Un chemin de $L_0$, est une famille de
plaques $(P_1,...,P_l)$ telle que $P_i\cap P_{i+1}$ n'est pas
vide, $l$ est appel\'e la longueur du chemin.

L'application de croissance $\gamma_{U_i}$ d\'efinie sur  ${\N}$,
est d\'efinie par $\gamma_{U_i}(n)$, est le nombre de plaques qui
peuvent \^etre jointes \`a $P_i$ par un chemin de longueur
inf\'erieur \`a $n$. La croissance de $L_0$ est celle de
l'application $\gamma_{U_i}$.

\medskip

Supposons que $M$ et la feuille $L_0$ soient  compactes, et soient
$T$ une transversale locale et connexe de ${\cal L}$ en $L_0$, et
$\hat L_0$, le rev\^etement unuversel de $L_0$. La m\'ethode de
suspension de Haefliger permet de repr\'esenter un voisinage de
$L_0$, comme quotient de $\hat L_0\times T$ par l'action $\rho$ du
groupe fondamental $\pi_1(L_0)$ d\'efinie par les transformations
du rev\^etement sur le premier facteur, et l'holonomie sur le
second.

\medskip

{\bf Proposition 2.1}

{\it Supposons que la feuille $L_0$ soit compacte, alors la
croissance de toute feuille $L_1$, appartenant au voisinage $V$ de
$L_0$, d\'efinie au paragraphe pr\'ec\'edent est major\'ee par
celle de $\pi_1(L_0)$.}

\medskip

{\bf Preuve.}

La transversale utilis\'ee pour construire la suspension de
Haefliger peut \^etre choisie compacte. On peut alors appliquer la
proposition $1.29$ de [6].

\medskip

{\bf Proposition 2.2.}

{\it Supposons que la croissance des feuilles de ${\cal L}$ dans
un voisinage de $L_0$ soit polynomiale, alors $\pi_1(L_0)$ est
polycyclique.}

\medskip

{\bf Preuve.}

Supposons que $\pi_1(L_0)$ ne soit pas polycyclique, alors
l'alternative de Tits, implique l'existence d'un sous-groupe libre
de $\pi_1(L_0)$ de rang sup\'erieur \`a $2$. On en d\'eduit que
l'holonomie infinit\'esimale de $L_0$ contient aussi un
sous-groupe libre de rang sup\'erieur \`a $2$. La croissance de
toute feuille du voisinage de $L_0$ ne peut \^etre polynomiale.
Une autre preuve de ce r\'esultat se trouve dans l'article de
Thurston et Plante.

\medskip

Dans une version ant\'erieure de ce travail, nous avons
conjectur\'e que toute vari\'et\'e affine compacte et compl\`ete
$(N,\nabla_N)$ est la feuille d'un feuilletage lagrangien, et les
feuilles d'un voisinage  de $N$, ont une croissance polynomiale.
William Goldman nous a communiqu\'e un contre exemple \`a cette
conjecture que voici:

Consid\'erons la suspension $M$ du tore de dimension $2$, $T^2$
munie de sa structure riemannienne plate standard au-dessus du
cercle $T^1$, par l'action de ${\Z}$ d\'efinie par un \'el\'ement
de $Gl(2,{\Z})$ hyperbolique, c'est \`a dire ayant des valeurs
propres distinctes en valeurs absolues de $1$.  Le fibr\'e
cotangent $T^*M$ de $M$ peut \^etre compactifi\'e par l'action du
tore $T^3$ agissant fibre \`a fibre par des translations. La
vari\'et\'e $M$ est donc feuille d'un feuilletage lagrangien. La
croissance de son groupe fondamental est exponentielle.

\medskip

{\bf Remarque.}

Il existe des vari\'et\'es affine compactes dont le groupe
fondamental n'est pas polycyclique. C'est le cas du produit du
cercle par une surface de genre sup\'erieur \`a $2$.

\medskip

{\bf Proposition 2.3.}

{\it Soit $(L_0,\nabla_{L_0})$ une vari\'et\'e affine compacte
dont l'holonomie lin\'eaire est incluse dans $Gl(n,{\Z})$, alors
$L_0$ est la feuille d'un feuilletage lagrangien d\'efini sur une
vari\'et\'e symplectique compacte.}

{\bf Preuve.}

Soit $\hat L_0$ le rev\^etement universel de $L_0$, le fibr\'e
cotangent $T^*L_0$ de $L_0$, est le quotient de $\hat L_0\times
{\R}^n$ par l'action de $\pi_1(L_0)$ d\'efinie par sur $\hat L_0$ par les transformations
du rev\^etement et sur ${\R}^n$ par:

$$
 h:\pi_1(L_0)\rightarrow Gl(n,{\R}),
$$

 $$
 \gamma\rightarrow
^t(L(h_{L_0})(\gamma^{-1})),
$$

o\`u $L(h_{L_0})$ est l'holonomie lin\'eaire de
$(L_0,\nabla_{L_0})$.

Soient $e_1$,..,$e_n$ une base de ${\R}^n$, et $t_{e_i}$ la
translation de direction $e_i, 1\leq i\leq n$. Le quotient de
$T^*L_0$ par $t_{e_1},..,t_{e_n}$ fibre \`a fibre est une
vari\'et\'e compacte $(N,\nabla_N)$, suspension du tore $T^n$
au-dessus de $L_0$. La structure symplectique du cotangent $T^*L_0$ se projette
en une structure symplectique de $N$.

\medskip

Soit $(M,\nabla_M)$ une vari\'et\'e symplectique compacte munie
d'un feuilletage lagrangien ${\cal L}$, dont $(L_0,\nabla_{L_0})$
est une feuille compacte et compl\`ete. Il existe un voisinage $U$
de $L_0$ symplectomorphorphe \`a la section nulle du fibr\'e
cotangent de $L_0$. L'holonomie infinit\'esimale de ${\cal L}$ et
celle de $\nabla_{L_0}$ sont mis en dualit\'e par la structure
symplectique. On en d\'eduit que la conjecture d'Auslander est
v\'erifi\'ee si la croissance combinatoire des feuilles du
feuilletage horizontal de $T^*L_0$ est polynomiale.

\medskip

Milnor a demand\'e si le groupe fondamental d'une vari\'et\'e
affine est polycyclique. Margulis a construit un exemple de
vari\'et\'e affine de dimension $3$ dont le groupe fondamental est
un groupe libre \`a deux g\'en\'erateurs. D'autres exemples de
vari\'et\'es affines dont le groupe fondamental est libre ont
\'et\'e construits par Goldman, Drumm et Charette. Ces
vari\'et\'es sont-elles des feuilles de feuilletages lagrangiens?

Le groupe fondamental d'une vari\'et\'e affine compacte non
compl\`ete peut contenir un sous-groupe libre \`a deux
g\'en\'erateurs. Ces vari\'et\'es sont-elles des feuilles de
feuilletages lagrangiens ?

\medskip

{\bf Examples de vari\'et\'es affines compactes feuilles de
feuilletages lagrangiens.}

\medskip

 On d\'efinit deux feuilletages lagrangiens sur le tore
de dimension $4$, $T^4$ ayant des feuilles lagrangiennes
compactes.

Le premier exemple est le tore de dimension $4$, $T^4$ muni de sa
structure euclidienne plate standard, il est le quotient de
${\R}^4$ par les translations $t_{e_1},...,t_{e_4}$, o\`u
$e_1,...,e_4$ est une base de ${\R}^4$. La forme symplectique est
la forme $\omega_0=dx_1\wedge dx_3 + dx_2\wedge dx_4$. Le feuilletage
lagrangien est l'image du feuilletage de ${\R}^4$ par espaces
affines parral\`eles \`a $vect(e_1,e_2)$ par l'application
rev\^etement.

\bigskip

Le second exemple est le quotient de ${\R}^4$ par les
transformations affines

$$
h_1(x_1,x_2,x_3,x_4)=(x_1+1,x_2,x_3,x_4)
$$

$$
h_2(x_1,x_2,x_3,x_4)= (x_1+x_2,x_2+1,x_3,-x_3+x_4)
$$

$$
h_3(x_1,x_2,x_3,x_4)= (x_1,x_2,x_3,x_4+1)
$$

$$
h_4(x_1,x_2,x_3,x_4)= (x_1+x_2,x_2,x_3+1,-x_3+x_4)
$$

La structure symplectique est la projection de $\omega_0$. Le
feuilletage lagrangien est l'image du feuilletage de ${\R}^4$ par
espaces affines parral\`eles \`a $vect(e_1,e_2)$ par l'application
rev\^etement. La feuille compacte est la projection de celle
passant par $0$.

\bigskip

Le quotient de ${\R}^4$ par les transformations affines $h_1$,
$h_2$ $t_3$ et $t_4$ est une vari\'et\'e symplectique. Le
feuilletage lagrangien est l'image du feuilletage de ${\R}^4$ par
espaces affines parral\`eles \`a $vect(e_1,e_2)$ par l'application
rev\^etement.

\bigskip

 Consid\'erons une vari\'et\'e symplectique $M$ de
dimension $2n, n\geq 2$. Donaldson a montr\'e qu'il existe des
sous-vari\'et\'es symplectiques $N_m$, de codimension $2m$ dans
$M$ dont le groupe fondamental est isomorphe \`a celui de $M$, si
la dimension de $N_m$ est sup\'erieure ou \'egale \`a $4$. Le
plongement naturel de $N_1\rightarrow M$ induit in \'epimorphisme
de groupes fondamentaux. Les vari\'et\'es $N_m$ peuvent \^etre
choisies disjointes de toutes sous-vari\'et\'es lagrangienne $N$, en contenant tout point contenu
dans $M-N$. Auroux a d\'emontr\'e que leur classe d'isomorphisme
est un invariant de la structure symplectique. On en d\'eduit:

\bigskip

{\bf Proposition 2.4.}

{\it Soit $(M,\nabla_M)$ une vari\'et\'e affine compacte et
compl\`ete dont l'holonomie lin\'eaire est contenue dans dans
$Gl(n,{\Z})$, alors le compactifi\'e du fibr\'e cotangent
$(N,\nabla_N)$ d\'efinie \`a la proposition pr\'ec\'edente \`a une
sous-vari\'et\'e symplectique de dimension $4$ dont le groupe
fondamental est isomorphe \`a $\pi_1(N)$ qui est une extension de
${\Z}^m$ par $\pi_1(M)$.}

\bigskip

{\bf Proposition 2.5.}

{\it Sous les hypoth\`eses de la proposition pr\'ec\'edente,
consid\'erons   $N'_m$ une composante connexe de $p^{-1}(N_m)$ par la
projection $p:{\R}^{2n}\rightarrow N$, alors $N'_m$ n'est pas
contractible.}

\bigskip

{\bf Preuve.}

Supposons que $N'_m$ soit contractible, l'application
$\pi_1(N_m)\rightarrow \pi_ 1(N)$ est un isomorphisme. Ceci
implique que $\pi_1(N)$ pr\'eserve $\hat N_m$ qui est
contractible. Les vari\'et\'es $N$ et $N_m$ \'etant des espaces
d'Eilenberg Mclane, on d\'eduit que la dimension de $N_m$ est
sup\'erieure ou \'egale \`a celle de $N$.

\bigskip

{\bf Remarque.}

Le morphisme $\pi_1(N_1)\rightarrow \pi_1(N)$ est un
\'epimorphisme strict car la dimension cohomologique de
$\pi_1(N_1)$ est $2$.

\bigskip

{\bf Proposition 2.6.}

{\it L'ensemble $p^{-1}(N_m)$ est connexe, si $m>1$ est n'est pas
alg\'ebrique.}

\bigskip

{\bf Preuve.}

L'ensemble $p^{-1}(N_m)$ est connexe car il est stable par
$\pi_1(N)$. Supposons qu'il soit alg\'ebrique, ceci entraine qu'il
est stable par l'adh\'erence de Zariski de $\pi_1(N)$ dans le
groupe des transformations affines $Aff({\R}^{2n})$ de ${\R}^{2n}$
qui agit transitivement sur ${\R}^{2n}$. On en d\'eduit une
contradiction.

\bigskip

{\bf D\'efinition 2.1.}

Soit $N$ une vari\'et\'e diff\'erentiable, une fibration de
Lefstchetz d\'efinie sur $N$, est une application $h:N\rightarrow
U$ de $N$ dans une surface de Riemann $U$ v\'erifiant les
propri\'et\'es suivantes:

 $h$ est une submersion sur $N-\{n_1,..,n_p\}$

 Sur une carte $(U_i,\phi_i)$ de $n_i$, $h$ est  une fonction de morse:
il existe un voisinage $(V_i,\psi_i)$ de $h(n_i)$ tel que

$$
\psi_i\circ h\circ \phi^{-1}(x_1,..,x_{n})=\sum {x_i}^2
$$

Un pinceau symplectique sur $M$, est d\'efini par un recouvrement
de $M$ par des vari\'et\'es symplectiques de codimension $2$,
telles que:

Il existe une vari\'et\'e de codimension $4$, $N$, telle que
l'intersection de deux vari\'et\'es distinctes de la famille
pr\'ec\'edentes est $N$,

L'\'eclat\'ee de $M$ le long de $N$ est  une fibration de
Leftschetz

\bigskip

Donaldson a montr\'e l'existence d'un pinceau symplectique sur
toute vari\'et\'e symplectique compacte. En dimension $4$ ce
pinceau est une fibration de Leftschetz
 dont l'image est la sph\`ere, les
fibres r\'eguli\`eres sont des surfaces diff\'eomorphes \`a $C_l$
la surface compacte de genre $l$, La fibre en $n_i$ est obtenu en
contractant un lacet $c_i$ de $C_l$, le groupe fondamental de $N$
est le quotient de $\pi_1(C_l)$ par le sous-groupe normal
engendr\'e par les  lacets.

\bigskip

{\bf Proposition 2.7.}

{\it Supposons que le pinceau sur une vari\'et\'e symplectique
quelconque $M$ soit une fibration, alors $M$ est diff\'eomorphe
\`a un produit $S^2\times N$.}

\bigskip

{\bf Preuve.}

L'horthogonal des fibres de la fibration symplectique est
int\'egrable. On conclut grace \`a un th\'eor\`eme de Ehresmann.

\bigskip

{\bf Proposition 2.8.}

{\it La restriction de la connection $\nabla_N$ d\'efinissant la
structure affine de $N$, induit sur chaque sous-vari\'et\'e $N_i$
une connection sym\'etrique $\nabla_{N_i}$.}

\bigskip

{\bf Preuve.}

Soient $x$ un \'el\'ement de $N_i$, et $X_1$, $X_2$ deux champs de
vecteurs de $N_i$. La restriction du fibr\'e tangent $TN$ de $N$
\`a $N_i$, est la somme $TN_i\oplus TN'_i$, o\`u $TN'_i$ d\'esigne
l'orthogonal de $TN_i$ par la forme symplectique. Consid\'erons un
atlas affine $(U_i)_{i\in I}$ de $(N,\nabla_N)$ pour tout champs
de vecteurs $Y_1$ et $Y_2$ de $N$, ${\nabla_N}_{{Y_1}_{\mid
U_i}}{Y_2}_{\mid U_i}=d{Y_2}_{\mid U_i}({Y_1}_{\mid U_i})$,
${\nabla_{N_i}}_{{X_1}_{\mid U_i}}{X_2}_{\mid U_i}$ est la
projection de $d{X_2}_{\mid U_i}({X_1}_{\mid U_i})$ sur $TN_i$
parall\`ement \`a $TN'_i$.

\bigskip

{\bf Remarque.}

La vari\'et\'e $N_i$ peut \^etre consid\'er\'ee comme une fibre
d'un pinceau symplectique de $N_{i+1}$.

\bigskip

{\bf Proposition 2.9.}

{\it La restriction de la forme symplectique $\omega_i$ de $N$ \`a
$N_i$, est parall\`ele relativement \`a $\nabla_{N_i}$.}

\bigskip

{\bf Preuve.}

Supposons que $\omega_{i+1}$ soit parall\`ele relativement \`a
$\nabla_{N_{i+1}}$. Consid\'erons une sous-vari\'et\'e de
codimension $4$, $L_i$ de $N_{i+1}$ telle que $N_i$ soit la
feuille d'un pinceau sur $N_{i+1}$ dont les feuilles contiennent
toutes $L_i$. La restriction du pinceau $P_i$ \`a $N_{i+1}-L_i$
est un feuilletage singulier ${\cal F}_i$.  Soient $X,Y,Z$ trois
champs de vecteurs de $N_i$, pour tout point $u$ de $N_i-L_i$, il
existe un voisinage $U$ de $u$ dans $N_{i+1}-L_i$ tel que la
restriction de ${\cal F}_i$ \`a $U$ soit simple, $U=F_0\times T$,
et les feuilles de la restriction du feuilletage \`a $U$ sont
$F_0\times x$. Soient trois champs de vecteurs tangents aux
feuilles $X'$,$Y'$ et $Z'$ dont les restrictions respectives \`a
${N_{i}-L_i}$ sont $X$, $Y$, $Z$. On suppose  que $X'=(X,0)$,
$Y'=(Y,0)$ et $Z'=(Z,0)$ sur $U$. La forme symplectique
$\omega_{i+1}$ est parall\`ele par rapport \`a $\nabla_{i+1}$ ceci
s'\'ecrit:

$$
X'.\omega_{i+1}(Y',Z')-\omega_{i+1}({\nabla_{N_{i+1}}}_{X'}Y',Z')-
\omega_{i+1}(Y',{\nabla_{N_{i+1}}}_{X'}Z')=0.
$$

De l'expression  des champs $X'$, $Y'$ et $Z'$  tangents aux
feuilles de ${\cal F}_i$, on d\'eduit que sur $N_i-L_i\cap U$, on
a:

$$
X.\omega_{i}(Y,Z)-\omega_{i}({\nabla_{N_{i}}}_{X}Y,Z)-
\omega_{i}(Y,{\nabla_{N_{i}}}_{X}Z)=0.
$$
 Cette relation est v\'erifi\'ee sur $N_i$.

 \bigskip

{\bf Reconstruction des vari\'et\'es affines symplectique \`a
l'aide des surfaces.}

\bigskip

Soit $TN'_1$, l'hortogonal du fibr\'e tangent $TN_1$ de $N_1$ pour
la forme symplectique. La connexion $\nabla_M$ induit sur $TN'_1$
une connexion $\nabla_{N'_1}$.

Un probl\`eme important est de reconstruire les vari\'et\'es
affines symplectiques compactes \`a partir d'une surface.

\medskip

 Toute vari\'et\'e
symplectique affine se construit ainsi:

On consid\`ere une surface $N_1$, un fibr\'e $TN'_1$ au-dessus de
$N_1$  telle que la somme des fibr\'es $V=TN'_1\oplus TN_1$ o\`u
$TN_1$ d\'esigne le fibr\'e tangent de $N_1$ soit plate et
symplectique, et la restriction de la forme symplectique $\omega$
\`a $TN_1$ n'est pas d\'eg\'en\'er\'ee. On munit $V$ d'une
connexion sans courbure $\nabla_V$. La dualit\'e symplectique
induit des connexions $\nabla_{N_1}$, et $\nabla_{N'_1}$ sur
$TN_1$ et $TN'_1$. On suppose que la restriction de  la forme
symplectique sur $N_1$ soit parall\`ele par rapport \`a la
connexion $\nabla_{N_1}$. Consid\'erons un voisinage $U$ de la
section nulle de $TN'_1$, la connexion $\nabla_V$ est d\'efinie
par une famille de formes $\omega_i$ d\'efinies sur un
recouvrement $L_i$ de $N'_1$, \`a ce recouvrement est associ\'e un
recouvrement $U_i$ de $U$ par des ouverts $L_i\times D_i$, o\`u
$D_i$ est un disque de la fibre type de $TN'_1$, on d\'efinit sur
$U$ la connexion

$$
{\nabla_U}_XY=dY(X) + \omega_i(X)Y
$$

Cette d\'efinition a un sens car les groupes structuraux du
fibr\'e tangent de $U$ et de $V$ sont identiques. On suppose que
la courbure et la torsion de la connexion $\nabla_V$ est nulle. Le
groupe fondamental de $U$ est $\pi_1(N_1)$, on en d\'eduit une
repr\'esentation d'holonomie

$$
h:\pi_1(N_1)\longrightarrow Aff({\R}^{2n})
$$

La d\'evellopante de la structure affine de $U$ est un
diff\'eomorphisme

$$
D:\hat U\longrightarrow {\R}^{2n}
$$

 On suppose que l'image de $h$ est l'holonomie d'une vari\'et\'e
affine compacte $(M,\nabla_M)$ munie d'une forme symplectique
parall\`ele. La conjecture de Markus affirme qu'une telle
structure est compl\`ete.

\bigskip
\bigskip
{\bf Espace des modules des vari\'et\'e affines symplectiques.}

\bigskip
\bigskip
 Ce paragraphe est consacr\'e \`a l'\'etudes des relations
entre les espaces de modules des fibr\'es plats au-dessus des
espaces de Riemann et des espaces des modules des vari\'et\'es
affines symplectiques.

D'apr\`es le paragraphe pr\'ec\'edent, une vari\'et\'e affine
symplectique compacte est d\'efinie par une surface de Riemann
$N_1$, un rev\^etement $\hat N_1$ plong\'e dans ${\R}^{2n}$, et
sous-vari\'et\'e symplectique pour la forme symplectique
parall\`ele $\omega$, une repr\'esentation

$$
h:\pi_1(N_1)\longrightarrow Aff({\R}^{2n})
$$

v\'erifiant:

(i) la partie lin\'eaire pr\'eserve la forme symplectique
$\omega$.

(ii) L'image $h(\pi_1(N_1))$ pr\'eserve $\hat N_1$ et le quotient
de ${\R}^{2n}$ par $h(\pi_1(N_1))$ est une vari\'et\'e affine
compacte.

\bigskip
\bigskip

{\bf Th\'eor\`eme 2.10.}

{\it Soit $(N,\nabla_N)$ une vari\'et\'e affine compacte dont le
groupe d'holonomie lin\'eaire est contenue dans $Gl(n,{\Z})$, le
groupe fondamental de $N$ est polycyclique. }

\bigskip

{\bf Preuve.}

\bigskip

On a construit une vari\'et\'e affine $(N',\nabla_{N'})$ qui est le
quotient du fibr\'e cotangent de $N$ par ${\Z}^n$, la vari\'et\'e
$N'$ est munie d'une forme symplectique $\omega_{N'}$ parral\`ele
relativement a la connexion affine. Il suffit de montrer que le
groupe fondamental de $N'$ est polycyclique.

Soit $N_1$ une surface de Donaldson contenue dans $N'$ telle que
l'application naturelle $\pi_1(N_1)\rightarrow \pi_1(N')$ soit
surjective.  Il a \'et\'e montr\'e par Weinstein que les classes de
Chern d'ordre impaire pour toutes structures pseudo-complexes
sous-jacentes a un fibr\'e symplectiques qui contient un
sous-fibr\'e lagrangien sont nulles. Il en r\'esulte que les classes
de Chern d'ordre impaires de $TN'$ sont nulles pour toutes
structures pseudo-complexes subordonn\'ee car les fibres de la
projection $N'\rightarrow N$ sont des tores lagrangiens. On en
d\'eduit que la premi\`ere classe de Chern de $TN'_{\mid N_1}$ la
restriction de $TN'$ \`a $N_1$ est nulle. Ceci implique que
$TN'_{\mid N_1}$ est un fibr\'e complexe trivial.

Les fibr\'e complexes  au-dessus de $N_1$ sont classifi\'es par
leurs premi\`ere classe de Chern. Le fibr\'e normal de $N_1$ peut
\^etre r\'ealis\'e comme un fibr\'e holomorphe. D'apr\`es le
th\'eor\`eme d'identication des sous-vari\'et\'e symplectiques \`a
fibr\'e normal isomorphes, on peut supposer que la structure presque
complexe sur un voisinage  $U$ de $N_1$ est int\'egrable.

Soient $v_1,..,v_n$ une trivialisation complexe de $TN'_{\mid N_1}$,
on peut \'etendre les vecteurs $v_1,..,v_n$ sur un voisinage de
$N_1$ en des champs de vecteurs complexes en suivant les
g\'eod\'esiques suivant l'orthogonal pour la structure symplectique
du fibr\'e tangent de $N'_1$.

   $N_1$ \'etant compacte, les
coordonn\'ees
 de $[v_i,v_j]$  dans la trivialisation complexe $v_l$,   sont des fonctions
 holomorphes de $N_1$ d'apr\`es le th\'eor\`eme du
 prolongement analytique, elle sont donc constantes.  On d\'efinit ainsi
 une alg\`ebre de Lie ${\cal H}$ et on note $H$ le groupe de Lie simplement connexe correspondant.

 Soit $N'_1$ une composante connexe de $p^{-1}(N_1)$, o\`u $p:\hat N\rightarrow N'$
 est la projection rev\^etement et $\hat N$ le rev\^etement universel de
 $N'$, on d\'enote par $U'$ une composante connexe de $p^{-1}(U)$
 contenant $N'_1$. Les champs $v_1,..,v_n$ se rel\`event sur $U'$
 en des champs holomorphes $v'_1,..,v'_n$ holomorphes.
 Le groupe $H$ est aussi d'alg\`ebre de Lie engendr\'ee par $v'_1,..,v'_n$ et agit sur
 $U'$
  comme un pseudogroupe
 transitif,
 Soit $x_0$ un
\'el\'ement
 de $U'$, et $\gamma$ un \'el\'ement de $\pi_1(N)$, il existe
  un \'el\'ement $h_{\gamma}$
 de $H$ tel que $\gamma(x_0)=h_{\gamma}(x_0)$. Le voisinage $U'$ s'identifie
 \`a un ouvert de $U"$ de $H$. Pour tout \'el\'ement $x$ de $U'$,
 il existe un \'el\'ement
 $n$ de $U'$ tel que $n(x_0)=x$, $nh_{\gamma}(x_0)=\gamma(n(x_0))$ puisque $\pi_1(M)$ pr\'eserve ${\cal H}$,
 on d\'eduit que l'action de $\pi_1(M)$ sur $N'_1$ coincide avec l'action \`a droite d'un
 sous-groupe discret de $H$ sur $U'$ isomorphe \`a $\pi_1(N)$.

 Consid\'erons la d\'ecomposition d'Iwasawa de $H=CAN$, o\`u $H$ est compact, $A$ ab\'elien
 et $N$ nilpotent. $C\cap \pi_1(N)$ est fini, car $\pi_1(N)$ est discret,  quitte \`a
 se placer sur un rev\^etement fini de $N$ tel que les valeurs propres de l'holonomie lin\'eaire
 n'ont pas d'ordre fini (voir Raghunatan), on peut supposer que $C\cap \pi_1(M)$ est
 l'identit\'e. Le quotient $H/C\pi_1(N)$ est une vari\'et\'e Eilenberg Maclane qui a le
 m\^eme type d'homotopie que la vari\'et\'e d'Eilenberg Maclane $N$, on d\'eduit que $H/C\pi_1(N)$ et $N$ ont m\^eme dimension,
 car $N$ est compacte, ceci implique $C$ est l'identit\'e et que $H$ est r\'esoluble,
 on en d\'eduit que $\pi_1(N)$ est r\'esoluble.

\bigskip

\bigskip

{\bf Corollaire.}

{\it Les conjectures d'Auslander et de Markus sont vraies si
l'holonomie lin\'eaire est contenue dans $Gl(n,{\Z})$.}

\bigskip
\bigskip

{\bf Genre de la surface de Donaldson.}

\medskip

\medskip

 Soit $(N,\nabla_N,\omega_N)$ une
vari\'et\'e affine compacte compl\`ete munie d'une forme
symplectique $\omega_N$, le but de cete partie est de d\'eterminer
le genre de la surface de Donaldson $N_1$ en fonction de la
dimension $2n$ de $N$

\medskip

{\bf Proposition.}

{\it  Le genre de la surface de Donaldson $N_1$ de
$(N,\nabla_N,\omega_N)$ est sup\'erieur \`a $n$, $dim(N)=2n$ si
$(N,\nabla_N)$ est compl\`ete et son groupe fondamental est
r\'esoluble.}

\medskip

{\bf Preuve.}

Sous les hypoth\`eses du th\'eor\`eme, $(N,\nabla_N)$ a un
rev\^etement fini $(N',\nabla_{N'})$ qui est le quotient d'un
groupe de Lie r\'esoluble par un r\'eseau $\Gamma$. Le nombre
minimal de g\'en\'erateur de $\Gamma$ est $2n$, La vari\'et\'e de
Donaldson $N_1$ se rel\`eve en une vari\'et\'e de Donaldson $N'_1$
de $N'$, comme l'application $\pi_1(N'_1)\rightarrow \pi_1(N)$ est
surjective, on en d\'eduit que le nombre de g\'en\'erateurs de
$N'_1$ donc de $N_1$ est sup\'erieur \`a $2n$, et par suite le
genre de $N_1$ est sup\'erieur \`a $n$.

 {\bf 3. Conjecture de Markus et feuilletage
transversallement mesurable.}

\bigskip
\bigskip

Markus a conjectur\'e qu'une vari\'et\'e affine compacte
$(M,\nabla_M)$ est compl\`ete si est seulement si elle est
unimodulaire. Ceci signifie que le groupe d'holonomie lin\'eaire
est un sous-groupe de $Sl(n,{\R})$.

Supposons que la vari\'et\'e affine compacte $(L_0,\nabla_{L_0})$
est la feuille d'un feuilletage lagrangien ${\cal L}$ d\'efini sur
la vari\'et\'e symplectique $(M,\omega)$. Consid\'erons un
voisinage $V$ de $L_0$ diff\'eomorphe au quotient de $\hat
L_0\times T$, o\`u $T$ d\'esigne une transversale lagrangienne, et
la restriction de ${\cal L}$ \`a $V$ est le feuilletage
suspension. L'action produit de $\pi_1(M)$ est d\'efini sur $\hat
L_0$ par les transformations de rev\^etements et sur $T$ par
l'holonomie du feuilletage. Puisque la dualit\'e symplectique
d\'efinie un isomorphisme entre l'holonomie infinit\'esimale de
${\cal L}$ et l'holonomie lin\'eaire de $(L_0,\nabla_{L_0})$, on
d\'eduit que l'holonomie de ${\cal L}$ pr\'eserve une m\'esure de
$T$. Ceci conduit \`a la conjecture suivante:

\bigskip

{\bf Conjecture 3.1.}

{\it Une vari\'et\'e affine  compacte $(L_0,\nabla_{L_0})$ est
compl\`ete si et seulement si elle est la feuille d'un feuilletage
lagrangien muni d'une mesure transverse dont le support contient
un voisinage de $L_0$.}

\bigskip

Plante a montr\'e que la croissance d'une feuille contenue dans le
support d'une m\'esure transverse  d'un feuilletage de codimension
$1$ est polynomiale. Ce r\'esultat n'est pas vrai en codimension
sup\'erieure comme le montre certaines suspensions de $S^2$.

\bigskip

Une vari\'et\'e affine compacte $(L_0,\nabla_{L_0})$ et compl\`ete
dont le groupe fondamental est polycyclique est unimodulaire. Une
vari\'et\'e affine compacte et compl\`ete est-elle la feuille d'un
feuilletage lagrangien d'une vari\'et\'e  symplectique compacte
munie d'une m\'esure transverse dont le support contient un
voisinage satur\'e de $L_0$, et telle que la croissance des
feuilles appartenant \`a ce voisinage soit polynomiale?

\bigskip

{\bf Proposition 3.2.}

{\it Soit $(M,\nabla_M,\omega_M)$ une vari\'et\'e affine compacte
muni d'un feuilletage lagrangien ${\cal L}$ ayant une feuille
compacte $L_0$ dont le groupe fondamental est un sous-groupe
normal de $\pi_1(M)$, alors $M$ est le quotient du fibr\'e
cotangent de $L_o$ par un groupe de symplectomorphismes affines.}

\bigskip

{\bf Preuve.}

Le quotient de ${\R}^{2n}$ par par la restriction de l'holonomie
\`a $\pi_1(L_0)$ est le fibr\'e cotangent de $L_0$, puisque
$\pi_1(L_0)$ est normal dans $\pi_1(M)$, on d\'eduit que
$(M,\nabla_M)$ est le quotient de $T^*L_0$ par
$\pi_1(M)/\pi_1(L_0)$.

\bigskip
\bigskip

{\bf Classification des vari\'et\'es symplectiques munies d'un
feuilletage lagrangiens.}

\bigskip
\bigskip

Soit $(M,\omega,{\cal L})$ une vari\'et\'e symplectique munie d'un
feuilletage lagrangien ${\cal L}$, ayant une feuille compacte
$(L_0,\nabla_{L_0})$. Supposons qu'il existe une action
parall\`ele du cercle $T^1$ sur $(L_0,\nabla_{L_0})$ pour la
connexion $\nabla_{L_0}$, Molino a \'etendu l'action du cercle en
une action hamiltonienne \`a $M$, et construit la r\'eduction de
Marsden-Weinstein $(M_1,\omega_1)$ pour cette action. Cette
derni\`ere vari\'et\'e est munie d'un feuilletage lagrangien
${\cal L}_1$ qui a une feuille symplectique compacte $L_1$, qui
est la base d'un fibr\'e affine d'espace total $L_0$.

Supposons que la feuille $L_0$ est diff\'eomorphe au tore de
dimension $n$, $T^n$, sa structure affine est d\'efinie par une
alg\`ebre associative et commutative ${\cal H}$, dont l'alg\`ebre
de Lie associ\'ee est celle du groupe commutatif $H$. Le tore
$T^n$ est le quotient de $H$ par un r\'eseau. Le groupe $H$ est
aussi la composante connexe du groupe des automorphismes affines
de $(L_0,\nabla_{L_0})$. Ceci permet d'identifier le produit
associatif de ${\cal H}$ \`a la restriction \`a ${\cal H}$ du
produit associatif de $aff({\R}^n)$ d\'efini par:

$$
(C,c).(D,d)=(CD,C(d))
$$

o\`u $C$, et $D$ d\'esignent des \'el\'ements de $gl(n,{\R})$ et
$c$, $d$ des \'el\'ements de ${\R}^n$.

Consid\'erons les g\'en\'erateurs $\gamma_1,...,\gamma_n$ de
$\pi_1(L_0)$, Il existe des \'el\'ements
$(C_1,c_1)$,...,$(C_n,c_n)$ de $aff({\R}^n)$ tels que
$\gamma_i=exp((C_i,c_i))$, on dira que la structure associative de
${\cal H}$ est rationnelle si le ${\Q}$ espace vectoriel
engendr\'e par $(C_1,c_1),...,(C_n,c_n)$ est stable par le produit
associatif.

\bigskip

{\bf Proposition 3.3.}

{\it Soit $(M,\nabla_M)$ une vari\'et\'e affine diff\'eomorphe au
tore, supposons que l'alg\`ebre associative qui d\'efinit sa
structure affine soit rationnelle, alors il existe une action du
cercle parall\`ele $S^1$ sur $(M,\nabla_M)$.}

\bigskip

{\bf Preuve.}

Soit $(C,c)$ un \'el\'ement du ${\Q}-$espace vectoriel engendr\'e
par $(C_1,c_1),...,(C_n,c_n)$ tel que $(C,c).(C,c)=(C^2,C(c))=0$.
$(C,c)=\alpha_1(C_1,c_1)+...+\alpha_n(C_n,c_n)$. Puisque les
$\alpha_i$ sont rationels, il existe un entier $p$ tel que
$p\alpha_i$ est un entier. L'\'el\'ement $exp(p(C,c))$ appartient
\`a $\pi_1(L_0)$, on conclut en utilisant Tsemo.

\bigskip

{\bf Corollaire 3.4.}

{\it Soit $(M,\omega)$ une vari\'et\'e symplectique munie d'un
feuilletage lagrangien ${\cal L}$ qui a une feuille compacte $L_0$
diff\'eomorphe au tore, et dont la structure affine est
compl\`ete, il existe une suite de vari\'et\'es affines
$(M_l,\omega_l),..(M_1,\omega_1)$ telles que $(M_l,\omega_l)$ est
$(M,\omega)$, est $(M_1,\omega_1)$ est un tore, de plus
$(M_i,\omega_i)$ est la r\'eduction de Marsden-Weinstein de
$(M_{i+1},\omega_{i+1})$ sur laquelle agit un cercle.}

\bigskip

{\bf Preuve.}

Il existe une translation $t_u$ appartenant \`a $\pi_1(L_0)$ car
la structure affine de $L_0$ est rationnelle. Cette action
s'\'etend en une action hamiltonienne de $(M,\omega)$. La
r\'eduction de Marsden-Weinstein de $(M,\omega)$ est la
vari\'et\'e symplectique $(M_1,\omega_1)$ munie d'un feuilletage
lagrangien ayant une feuille compacte $L_1$ dont la structure
affine est aussi rationnelle. On peut r\'eiterer le proc\'eder.

\bigskip

Il existe des structures affines sur le tore de dimension $2$ dont
le groupe d'holonomie ne contient pas de translations.

Consid\'erons le groupe de Lie de dimension deux $H$ de
$Aff({\R}^2)$ dont les \'el\'ements sont
$f_{s,t}(x,y)=(x+sy+{s^2\over 2}+t,y+s)$. Ce groupe agit
simplement et transitivement sur ${\R}^2$. Soit $h$ un r\'eel non
rationnel. Le sous-groupe $I$ de $H$ engendr\'e par $f_{h,0}$ et
$f_{1,1}$ est un r\'eseau de $H$. Le quotient de $H$ par $I$
d\'efinit une structure affine sur le tore dont l'holonomie ne
contient pas de translation.

\bigskip

Plus g\'en\'eralement consid\'erons une suite de vari\'et\'es
affines compactes $(L_n,\nabla_{L_n})\rightarrow...\rightarrow
(M_1,\nabla_{L_1})$ telle que l'application
$(L_{i+1},\nabla_{L_{i+1}})\rightarrow (L_i,\nabla_{L_i})$ soit
une fibration affine dont les fibres sont d\'efinies par une
action parall\`ele du cercle sur $(L_{i+1},\nabla_{L_{i+1}})$. On
suppose aussi que $(L_{i+1},\nabla_{L_{i+1}})$  est la feuille
d'un feuilletage lagrangien sur une vari\'et\'e symplectique
$(M_{i+1},\omega_{i+1})$. La r\'eduction de Marsden-Weinstein
$(M_i,\omega_i)$ pour l'extension de l'action du cercle est une
vari\'et\'e symplectique $(M_i,\omega_i)$ munie d'un feuilletage
lagrangien ayant $(L_i,\nabla_{L_i})$ comme feuille. Dans la
partie suivante on classifiera les germes des feuilletages au
voisinage de $L_i$.

\bigskip
\bigskip
{\bf Classification du voisinage d'une feuille compacte d'un
feuilletage lagrangien.}

\bigskip
\bigskip
Dans cette partie on rappelle la classification effectu\'ee par
Molino et Curras-Bosch qu'on g\'en\'eralise.

Consid\'erons une vari\'et\'e affine compacte et compl\`ete
$(L_0,\nabla_{L_0})$, $x_0$ un \'el\'ement de $L_0$, et
$h_{\nabla_0}$ l'holonomie de $(L_0,\nabla_{L_0})$, on identifie
$T_{x_0}L_0$ \`a $T^*_0{\R}^n$ muni de sa connexion plate
canonique. La partie lin\'eaire $h'_{\nabla_0}$ de l'holonomie
induit une r\'epr\'esentation $h'_{\nabla_0}:\pi_1(L_0)\rightarrow
T^* {\R}^n.$

Supposons aussi d\'efinie une repr\'esentation
$h_{x_0}:\pi_1(L_0)\rightarrow Diff({\R}^n)$. Celle-ci est
l'holonomie d'un feuilleutage lagrangien qui a pour feuille
compacte la vari\'et\'e affine $(L_0,\nabla_{L_0})$ si et
seulement si, pour tout \'el\'ement de $D_0$, l'ensemble des
fonctions diff\'erentiables d\'efinies sur un voisinage de $0$
telles que $f(0)=0$, on a

$$
h_{\nabla_{L_0}}(\gamma)\circ d_0=d_0\circ h_{x_0}(\gamma^{-1})^*,
\leqno(1)
$$

o\`u $h_{x_0}^*(\gamma)(f)=f\circ h_{x_0}(\gamma)$.

La repr\'esentation $h'_{\nabla_{L_0}}$ munie $T^*_0{\R}^n$ d'une
action de $\pi_1(L_0)$, on notera $H^*(\pi_1(L_0),T^*{\R}^n)$ les
espaces de cohomologies correspondant. L'obstruction radiante
$\gamma\rightarrow h_{\nabla_{L_0}}(0)$ d\'efinit un \'el\'ement
$[h'_{\nabla_{L_0}}]$ de $H^1(\pi_1(L_0),T^*{\R}^n)$.

La r\'epr\'esentation $h_{x_0}$ muni $D_0$ d'une structure de
module d\'efini par

$$
\gamma\circ f=f\circ h_{x_0}(\gamma^{-1})
$$

L'application $d_0$ induit un morphisme

$$
d_0^*: H^1(\pi_1(L_0),D_0)\rightarrow H^1(\pi_1(L_0),T^*_0{\R}^n)
$$

On en d\'eduit le th\'eor\`eme:

\bigskip

{\bf Th\'eor\`eme 3.5.}

{\it Les germes de feuilletages lagrangiens ayant
$(L_0,\nabla_{L_0})$ comme feuille compacte et pour holonomie
$h_{x_0}$ sont classifi\'es a symplectomorphismes pr\`es par les
\'el\'ements de ${d^*_0}^{-1}([h_{\nabla_0}]$.}

\bigskip

{\bf Classification des germes de feuilletages lagrangiens autour
d'une feuille compacte $(L_0,\nabla_{L_0})$ dont la structure
affine n'est pas forc\'ement compl\`ete.}

\bigskip

Dans cette partie on g\'en\'eralise la classification de Molino et
Curras-Bosch sans supposer que la structure affine
$(L_0,\nabla_{L_0})$ soit compl\`ete. Une autre classification a
\'et\'e obtenue par les auteurs pr\'ec\'edents.

\bigskip

On peut identifier un voisinage de $L_0$ dans $M$ \`a un voisinage
de la section nulle de $T^*L_0$, et \'etendre l'application
d\'eveloppement $D_0: \hat L_0\rightarrow {\R}^n$ en une
application $D:T^*\hat L_0\rightarrow {\R}^{2n}$ telle que pour
tout \'el\'ement $(x,y)$ de $T^*\hat L_0$, on a
$D(x,y)=(D_0(x),y)$. Cette identification est possible car
$T^*\hat L_0$ est un fibr\'e trivial. Consid\'erons une
transversale $T$ qu'on identifie \`a un ouvert contractible de la
fibre de $x_0$. Le feuilletage lagrangien dans $U_0$ est d\'efini
par une suspension. Son relev\'e sur  $U_0=L_0\times T$ est le
feuilletage dont les feuilles sont les sous-vari\'et\'es
$L_0\times y$. Pour chaque \'el\'ement $x$ de $L_0$, la dualit\'e
symplectique identifie $T^*(x\times T)$ \`a $T^*_x\hat L_0$. On
peut aussi identifier $T^*(x\times T)$ \`a $T^*(x_0\times T)$.

Soit $U_x$ un voisinage de $0$ dans $T_x L_0$ telle que la
restriction de $exp_{x_{\nabla_{L_0}}}=exp_x$, l'exponentielle
associ\'ee \`a la connexion affine \`a $U_x$ soit injective. La
dualit\'e symplectique permet d'identifier $U_x$ \`a un ouvert
$V_x$ de $T^*(x\times T)=T^*(x_0\times T)$. On d\'efinit ainsi une
carte $V_x\rightarrow \hat L_0$. En utilisant la m\'ethode bien
connue de construction d'une d\'eveloppante, on obtient une
application $D_0:\hat L_0\rightarrow T^*T_{x_0}$.

 On a

$$
   h_{\nabla_{L_0}}(\gamma)\circ d_{x_0}=d_{x_0}\circ h_{x_0}(\gamma^{-1})^*,
     \leqno(1)
$$

o\`u $h_{x_0}^*(\gamma)(f)=f\circ h_{x_0}(\gamma)$.

\bigskip

R\'eciproquement consid\'erons une vari\'et\'e affine compacte
$(L_0,\nabla_{L_0})$ et $h_{x_0}:\pi_1(L_0)\rightarrow
Diff({\R}^n)$ une repr\'esentation. D\'efinissons sur ${\R}^{2n}$
la forme symplectique dont le relev\'e par la d\'eveloppante est
la forme symplectique de $T^*\hat L_0$ relev\'ee de celle de
$T^*L_0$ par l'application rev\^etement. Comme au-dessus, la
dualit\'e symplectique permet d'identifier la d\'evelopante de la
structure affine $(L_0,\nabla_{L_0})$ \`a une application $\hat
L_0\rightarrow T^*{\R}^n$. Le fait que l'holonomie de
$(L_0,\nabla_{L_0})$ pr\'eserve le relev\'e de la structure
symplectique ${\R}^n\times \hat L_0$ signifie que:

$$
h_{\nabla_{L_0}}(\gamma)\circ d_0=d_0\circ h_{x_0}(\gamma^{-1})^*,
\leqno(1)
$$

On en d\'eduit une application:

 $$
d_0^*: H^1(\pi_1(L_0),D_0)\rightarrow H^1(\pi_1(L_0),T^*_0{\R}^n)
$$

On a le th\'eor\`eme de classification:

\bigskip

{\bf Th\'eor\`eme.}

{\it Les germes de feuilletages lagrangiens autour de
$(L_0,\nabla_0)$ dont l'holonomie du feuilletage est $h_{x_0}$
sont classifi\'es \`a symplectomorphismes pr\`es par les
\'el\'ements de ${d_0^*}^{-1}([h_{\nabla_{L_0}}]$.}

\bigskip

Revenons \`a la classification initiale, on supposera que la
vari\'et\'e affine $(L_i,\nabla_{L_i})$ est compacte et
compl\`ete. Etant donn\'e une classe d'isomorphisme de feuilletage
$e_1$ autour de $(L_1,\nabla_{L_1})$, on classifie les germes de
feuilletages autour de $(L_2,\nabla_{L_2})$ dont la r\'eduction de
Marsden-Weinstein est le germe $e_1$.

\bigskip

Soit $x_2$ un \'el\'ement de $L_2$ dont l'image par la projection
$L_2\rightarrow L_1$ est $x_1$. Un germe de feuilletage autour de
$(L_2,\nabla_{L_2})$ est d\'efinie par une repr\'esentation
$h_{x_2}:\pi_1(L_2)\rightarrow Diff({\R}^{n+1})$ v\'erifiant la
relation $(1)$. Supposons que la r\'eduction de Marsden-Weinstein
de ce germe est la vari\'et\'e $(M_1,\omega_1)$ munie un
feuilletage lagrangien ayant pour feuille $(L_1,\nabla_{L_1})$ et
tel que le germe autour de $(L_1,\nabla_{L_1})$ est $e_1$. On a la
proposition suivante:

\bigskip

{\bf Proposition 3.7.}

{\it Le diagramme suivant est commutatif:

$$
\matrix {\pi_1(L_2)&{\buildrel{h_{x_2}}\over{\longrightarrow}}&
Diff({\R}^{n+1})\cr\downarrow p_2 &\ \ \ \ & \downarrow \cr
\pi_1(L_1) &{\buildrel{h_{x_1}}\over{\longrightarrow}}&
Diff({\R}^n)}
$$
}
\bigskip

{\bf Preuve.}

Consid\'erons l'application moment $J:U_2\rightarrow {\R}$, on a
suppos\'e que $J^{-1}(0)$ contient $L_2$. La transversale $T_2$ au
feuilletage est diff\'eomorphe \`a un ouvert de ${\R}^{l_1+1}$. Il
suffit de remarquer qu'on peut supposer que l'application moment
est une fonction coordonn\'ee au voisinage de $0$. Pour le
prouver, on remarque que la fonction utilis\'ee pour \'etendre
l'action du cercle dans un voisinage de $L_2$ Voir Molino
Th\'eor\`eme 3.2 p. 186 doit d\'ependre uniquement d'une seule
coordonn\'ee dans un voisinage de $0$, puisque la fonction basique
associ\'ee coincide avec l'application moment dans un voisinage.

\bigskip

On va noter $e_2$ la classe d'isomorphisme de ce germe, et on va
munir l'ensemble des couples $(e_1,e_2)$, $e_1$ est fix\'e, d'une
structure de gerbe.

\bigskip

Consid\'erons le site $Et_{L_0}$ dont les objets sont les
rev\^etements de $L_0$ et les morphismes, les morphismes de
rev\^etements. A tout objet $e$ de $Et_{L_0}$ on associe la
cat\'egorie $symp(e,L_0)$ de germes de feuilletages lagrangiens
qui ont $e$ comme feuille, et telle que l'holonomie du feuilletage
d'un objet de $symp(e,L_0)$ en un point $x_e$ au-dessus de $x_0$
est donn\'ee par la restriction de $h_{x_0}$ \`a $\pi_1(e)$. Les
automorphismes des objets de $symp(e,L_0)$ sont les exponentiels
des champs hamiltoniens des fonctions basiques du feuilletage
lagrangien. Il est facile de montrer que la correspondance
$e\rightarrow symp(e,L_0)$ est une gerbe sur $Et_{L_0}$. Ceci
signifie que les axiomes suivanrs sont v\'erifi\'es:

- Pour toute fl\`eche $U\rightarrow V$ entre \'el\'ements de
$symp(e,L_0)$, on a un morphisme $r_{U,V}:symp(V,L_0)\rightarrow
symp(U,L_0)$ telle que $r_{U,V}\circ r_{V,W}=r_{U,W}$.

- Condition de recollement pour les objets.

Consid\'erons une famille couvrante $(U_i)_{i\in I}$ d'un objet
$U$ de $Et_{L_0}$, pour chaque $i$, un objet $x_i$ de
$symp(U_i,L_0)$. Supposons qu'il existe une application
$g_{ij}:r_{U_i\cap U_i,U_j}(x_j)\rightarrow r_{U_i\cap
U_i,U_i}(x_i)$ telle que $g_{ij}g_{jk}=g_{ik}$, alors il existe un
objet $x$ de $symp(e,L_0)$ tel que $r_{U_i,U}(x)=x_i$.

-Condition de recollement pour les fl\`eches.

Consid\'erons deux objets $P$ et $Q$ de $symp(L_0,L_0)$
l'application $U\rightarrow Hom(r_{U,L_0}(P),r_{U,L_0}(Q))$ est un
faisceau.

- Il existe une famille couvrante $(U_i)_{i\in I}$ de $Et_{L_0}$
telle que pour tout $i$ la cat\'egorie $symp(U_i,L_0)$ n'est pas
vide.

- Soit $U$ un objet de $Et_{L_0}$, pour tous objets $x$ et $y$ de
$symp(U,L_0)$, il existe une famille couvrante $(U_i)_{i\in I}$ de
$U$ telle que les objets $r_{U_i,U}(x)$ et $r_{U_i,U}(y)$ sont
isomorphes.

Toute fl\`eche de $symp(U,L_0)$ est inversible, et il existe un
faisceau en groupes commutatifs $H$ sur $Et_{L_0}$ tel que pour
tout objet de $symp(U,L_0)$, $Hom(x,x)=H(U)$, et cette
identification  commute avec les restrictions.

\bigskip

{\bf Remarque.}

La gerbe que nous venons d'\'etudier est une gerbe triviale comme
le montre le th\'eor\`eme de classification. Nous l'appelerons
$symp(L_0)$. Soit $u$ une section globale de $symp(L_0)$,
${L_0}_u$ le relev\'e de $L_0$ \`a $u$ \`a savoir la feuille de
$u$ qui se projette sur $L_0$. Consid\'erons une transversale
lagrangienne $T$ en ${L_0}_u$, la dualit\'e symplectique permet
d'identifier $f$ au feuilletage vertical de $T^*T_{x_0}$. Les
automorphismes de cet objet sont les automorphismes verticaux du
feuilletage. C'est un groupe commutatif puisque $T$ est
lagrangien, ces \'el\'ements sont les exponentiels de champs
hamiltoniens.

\bigskip
\bigskip
{\bf Le cas g\'en\'eral.}

\bigskip
\bigskip

Consid\'erons  $(L_l,\nabla_{L_l})\rightarrow...\rightarrow
(L_1,\nabla_{L_1})$ une suite de vari\'et\'es affine compl\`etes
telles que $f_i:(L_{i+1},\nabla_{L_{i+1}})\rightarrow
(L_i,\nabla_{L_i})$ est une fibration affine dont les feuille sont
les orbites d'une action parall\`ele du cercle. Le but de cette
partie est de classifier les suites de germes de feuilletages
lagrangiens $U_i$ telles que: $(L_{i+1},\nabla_{L_{i+1}})$ est une
feuille de $U_{i+1}$, est $L_i$ est une feuille la r\'eduction de
Marsden-Weinstein de $U_{i+1}$.

\bigskip

Le germe de feuilletage lagrangien de $U_i$  est d\'efini par une
repr\'esentation

$$
\pi_1(L_i)\longrightarrow Diff({\R}^n)
$$

v\'erifiant la relation $(1)$ tel que le diagramme suivant soit
commutatif

$$
\matrix
{\pi_1(L_{i+1})&{\buildrel{h_{x_{i+1}}}\over{\longrightarrow}}&
Diff'({\R}^{n+1})\cr\downarrow p_i &\ \ \ \ & \downarrow \cr
\pi_1(L_i)
&{\buildrel{h_{x_i}}\over{\longrightarrow}}&Diff({\R}^n)}
$$

O\`u $Diff'({\R}^{n+1})$ d\'esigne les automorphismes se
projettant sur ${\R}^n$. Si $e_1$ un objet de $symp(L_1)$, on
associe \`a $e_1$ la gerbe $symp(e_1,L_2)$ d\'efinie sur le site
$Et_{L_2}$, telle que pour tout objet $e_2$ de $Et_{L_2}$, les
objets de la cat\'egorie $symp(e_1,L_2)(e_2)$ sont les germes
$U_2$ de feuilletages lagrangiens qui ont $(L_2,\nabla_{L_2})$
comme feuille, et telle que la r\'eduction de Marsden-Weinstein
par l'action du cercle est un feuilletage lagrangien qui a $L_1$
comme feuille, et le germe autour de $L_1$ est isomorphe $e_1$.

\bigskip

Supposons d\'efinie la gerbe $symp(L_p,e_1,...,e_{p-1})$ et soit
$e_p$ un objet de cette gerbe. On va d\'efinir la gerbe
$symp(L_{p+1},e_1,..,e_p)$ au-dessus du site $Et_{L_{p+1}}$, dont
les objets sont les germes de feuilletages lagrangiens $U_{p+1}$
qui ont $L_{p+1}$ comme feuille, et tels que la r\'eduction de
Marsden-Weinstein par l'action du cercle a $L_p$ comme feuille et
le germe de ce feuilletage lagrangien autour de $L_p$ est
isomorphe \`a $e_p$.

\bigskip

Les gerbes $symp(L_{p+1},e_1,...,e_p)$ sont des gerbes triviales.
Pour tout objets $e_p$ et $e'_p$ de $symp(L_p,e_1,...,e_{p-1})$ et
tout morphisme $f:e_p\rightarrow e'_p$ il existe un foncteur
$f^*:symp(L_{p+1},e_1,...,e'_{p})\rightarrow
symp(L_{p+1},e_1,...,e_{p})$ tel qu'il existe un isomorphisme

$$
c(f,g): (fg)^*\rightarrow g^*f^*
$$

satisfaisant la condition de $1-$descente

$$
(Id*c(f,g)\circ c(fg,h) = c(g,h)*Id \circ c(f,gh)
$$

Supposons que la dimension de $L_1$ est $l_1$, la dimension de
$L_i$ est $l_1+i-1$. On va noter par $D_i$ les germes en $0$ des
fonctions diff\'erentiables de ${\R}^{l_1+i-1}$ qui s'annulle \`a
l'origine. L'holonomie en $x_i$ munie $D_i$ d'une structure de
$\pi_1(L_i)$ module en posant:

$$
\gamma\circ f=f\circ h_{x_i}(\gamma^{-1})
$$

La diff\'erentielle $d_i$ en l'origine induit un isomorphisme:

$$
d_i^*:H^1(\pi_1(L_i),D_i)\longrightarrow
H^1(\pi_1(L_i),T^*_0{\R}^{l_1+i-1})
$$

Les classes d'isomorphismes de feuilletages dont l'holonomie est
$h_{x_i}$ sont classifi\'ees par ${d_i^*}^{-1}([h_{\nabla_i}])$.

Etant donn\'e une section globale $e_i$ de $symp(L_i)$
classifi\'ee par la classe $c_i\in H^1(\pi_1(L_i),D_i)$, on va
d\'eterminer le cocycle classifiant des sections globales de
$symp(L_{i+1})$ qui donnent lieu \`a $e_i$.

La surjection $\pi_1(L_{i+1})\rightarrow \pi_1(L_i)$ et le
diagramme commutatif $(1)$ induisent le diagramme suivant:

$$
\matrix{H^1(\pi_1(L_{i+1}),D_{i+1})&\longrightarrow&
H^1(\pi_1(L_i),{\R}^{l_1+i})\cr f_i\downarrow && g_i\downarrow \cr
H^1(\pi_1(L_i),D_i)&\longrightarrow&
H^1(\pi_1(L_i),{\R}^{l_1+i-1})}.
$$

L'image des cocycles classifiant des \'el\'ements de
$symp(L_{i+1})$ par $f_i$ est $c_i$.

Soient $\gamma_1,...,\gamma_{l_1+i-1}$ les g\'en\'erateur de
$\pi_1(L_i)$, on supposera que si $j\leq i$,
$\gamma_1,...,\gamma_{l_1+j-1}$ sont les g\'en\'erateurs de
$\pi_1(L_j)$. Remarquons puisque nous avons suppos\'e que
$(L_2,\nabla_{L_2})$ est fibr\'e en cercle au-dessus de
$(L_1,\nabla_{L_1})$, on peut supposer que
$h_{x_2}(\gamma_{l_1+1})$ est l'identit\'e, c'est \`a dire est le
g\'en\'erateur de $\pi_1(L_2)$ qui pr\'eserve la fibre du fibr\'e
en cerle.

Soit $c_1$ le cocycle qui d\'efinit la classe d'isomorphisme du
germe $U_1$, $c_1(\gamma_1)$ est le germe d'une fonction
diff\'erentiable d\'efinit sur ${\R}^{l_1}$. Soit $\gamma'$ un
\'el\'ement de $\pi_1(L_2)$ au-dessus de $\gamma_1$,
$c_2(\gamma')$ est le germe d'une fonction diff\'erentiable de
${\R}^{l_1+1}$ au-desssus de $c_1(\gamma_1)$.

\bigskip

Consid\'erons $L(h_i)$, l'holonomie lin\'eaire de la vari\'et\'e
affine $(L_i,\nabla_{L_i})$, et \'ecrivons
${\R}^{l_1+1}={\R}^{l_1}\oplus {\R}$. Pour tout \'el\'ement
$\gamma\in \pi_1(L_2)$, $L(h_2(\gamma))$ d\'epend de $p_2(\gamma)$
et de la projection de $L(h_2(\gamma))$ sur ${\R}e_{l_1+1}$
parall\`element \`a ${\R}^{l_1}$ qui est un cocycle $d_2$ qui est
un cocycle pour l'action triviale de $\pi_1(L_2)$ sur
${\R}e_{l_1+1}$.

Nous allons supposer que les holonomies sont lin\'earisables dans
des voisinage des feuilles $L_i$. Nous allons d\'eterminer les
relations entre les cocycles qui d\'efinissent les \'el\'ements de
$symp(L_2,e_1)$ et le cocycle classifiant de $e_1$.

Pour tout \'el\'ement $e_2$ de $symp(e_1,L_2)$, nous allons
d\'eterminer les relations entre son cocycle classifiant et celui
de $e_1$.

Pour tout \'el\'ement $\gamma$ de $\pi_1(L_1)$, on consid\`ere une
fonction $c_2(\gamma)$ de ${\R}^{l_1+1}$ au-dessus de
$c_1(\gamma)$ ceci signifie qu'il existe une fonction
$f_2(\gamma)$ telle que

$$
c_2(\gamma)=c_1(\gamma)(x_1,..,x_{l_1})+f_2(\gamma)(x_{l_{1+1}})
$$

et pour $\gamma_{l_1+1}$, on consid\`ere un \'el\'ement
$c_2(\gamma_{l_1+1})$ qui se projette sur $0$, c'est \`a dire une
fonction de $x_{l_1+1}$ la $l_1+1-$cooedonn\'ee. Le fait que $c_2$
soit un cocycle implique que:

$$
c_2(\gamma_{l_1+1}\gamma)=c_2(\gamma_{l_1+1})+c_2(\gamma)
$$

Pour tout \'el\'ement $\gamma$ de $\pi_1(L_2)$,

$$
c_2(\gamma\gamma')=c_2(\gamma)+\gamma c_2(\gamma')=
$$

$$
c_1(p_2(\gamma))+f_2(\gamma)+ p_2(\gamma)c_1(p_2(\gamma)) +
f_2(\gamma')\circ d_2(\gamma)
$$

Ceci implique que

$$
f_2(\gamma\gamma')=f_2(\gamma)+f_2(\gamma')\circ d_2(\gamma)
$$

Plus g\'en\'eralement, le cocycle $c_i(h)$ d'un \'el\'ement de
$symp(L_i)$ est d\'eterminer par  un germe de fonction
diff\'erentiable $f_i(h)$ de ${\R}e_{l_1+i-1}$ tel que

$$
c_i(h)(\gamma)= c_{i-1}(p_i(\gamma)+f_i(h)(\gamma)
$$

qui satisfait

$$
 f_i(h)(\gamma\gamma')= f_i(h)(\gamma)+
 f_i(h)(\gamma)\circ d_i(\gamma).
$$

\bigskip
\bigskip
{\bf Bibliography.}

\bigskip
\bigskip

[1] Auroux, D. Techniques approximativement holomorphes et
invariants de monodromie en topologie symplectiques. Th\`ese
d'habilitation 2003.

\bigskip

[2] Charette, V. Goldman, W.

  Affine schottky groups and Crooked tilings. Proceedings of Crystallographic
  groups and their generalizations II. Kortrijk 1999.

\bigskip

[3] Curras-Bosch, C. Molino, P.

Holonomie suspension et classification pour les feuilletages
lagrangiens. C.R.A.S (326) 1317-1320

\bigskip

[4] Dazord, P.

Sur la geometrie des sous-fibres et des feuilletages lagrangiens.
Annales scient. E. Norm. Sup. Paris.

\bigskip

[5] Donaldson, Symplectic submanifolds and almost-complex
geometry. J. Differential Geometry (44) 1996 666-705.

\bigskip

[6] Donaldson, Lefstchetz pencils on symplectic manifolds.

J. Diffferential Geom. (53) 1999 205-236
\bigskip

[7] Fintushel, R. Stern, R. Symplectic surfaces in a fixed
homology class. J. Differential Geom. (52) 1999 203-222.

\bigskip

[8] Fried, D. Goldman, W.

Three-dimensional affine
 crystallographic groups. Advances in Math. 47 (1983), 1-49.

\bigskip

[9] Fried, D. Goldman, W. Hirsch, M.

 Affine manifolds with
 nilpotent holonomy. Comment. Math. Helv. 56 (1981) 487-523.

\bigskip

[10] Goldman, W. Hirsch, M. The radiance obstruction and parallel
forms on affine manifolds. Trans. Amer. Math. Soc. 286 (1984),
629-949. [6] Godbillon, C.

\bigskip

[11] Goldman, W. Hirsch, M. Affine manifolds and orbits of
Feuilletages. Etudes g\'eom\'etriques. Progress algebraic groups.
Trans. Amer. Math. Soc. 295 (1986), 175-198. in Mathematics, 98.

\bigskip

[12] Goldman, W. Geometric structure on manifolds and varieties of
representations. 169-198, Contemp. Math., 74

\bigskip

[13] Goldman, W. Private communication.

\bigskip

[14] Gunning, Lectures on vector bundles over surfaces.

\bigskip

[15] Koszul, J-L. Vari\'et\'es localement plates et convexit\'e.
Osaka J. Math. (1965), 285-290.

\bigskip

 [16] Koszul, J-L. D\'eformation des connexions localement plates.
Ann. Inst. Fourier 18 (1968), 103-114.

\bigskip

[17] Liberman, P. Sur les structures presques complexes et autres
structures infinit\'esiment r\'eguli\`eres.

Bulletin Soci\'et\'e Mat. France (83) 1955 195-224

\bigskip
 [18] Margulis, G.

 Complete affine locally flat manifolds with a
 free fundamental group. J. Soviet. Math. 134 (1987), 129-134.

\bigskip
 [19] Milnor, J. W.

 On fundamental groups of complete affinely flat
 manifolds, Advances in Math. 25 (1977) 178-187.

\bigskip
[20] Molino, P Lagrangian holonomy.

Analysis and geometry in foliated manifolds. (Santiago de
Compostela) 183-194.

\bigskip

[21] Plante, J.F. Polycyclic groups and transversely affine
foliations.

J. Diff. Geometry 35 (1992) 521-534.

\bigskip

[22] Plante, J. F. Foliations with measures preserves holonomy.

Ann. of Math. 102 (1975) 327-361.

\bigskip

[23] Plante J.F Thurston, W. Polynomial growth in holonomy
foliations.

Comment. Math. Helv. 51 (1976) 567-584.

\bigskip

[24] Sullivan, D. Thurston, W. Manifolds with canonical
 coordinate charts: some examples. Enseign. Math 29 (1983), 15-25.

\bigskip

 [25] Tsemo, A. Dynamique des vari\'et\'es affines. J. London Math.
 Soc. 63 (2001) 469-487.

\bigskip

[26] Weinstein, A. Symplectic manifolds and their lagrangian
submanifolds.

Advances in Math 6  329-346.

\end{document}